\documentclass{amsart}
\newcommand{\note}{\noindent {\bf Notation. }}
\newcommand{\remark}{\noindent {\bf Remark. }}

\newcommand{\ws}{\hspace{4pt}}
\newtheorem{theorem}{Theorem}

\newtheorem{proposition}{Proposition}
\newtheorem{lemma}{Lemma}

\usepackage[]{amssymb,amsopn}
\usepackage[]{amsmath}
\usepackage[]{eucal}
\usepackage[]{latexsym}
\usepackage[]{mathtools}

\begin{document}

\title[Asymptotic Recurrence Coefficients]{Asymptotics for Recurrence Coefficients of X1-Jacobi Polynomials and Christoffel Function}
\author{\'A. P. Horv\'ath }

\subjclass[2010]{33C47,33C45}
\keywords{exceptional orthogonal polynomials, recurrence relation, Chistoffel function}
\thanks{Supported by the NKFIH-OTKA Grant K128922.}

\begin{abstract} Computing asymptotics of the recurrence coefficients of X1-Jacobi polynomials we investigate the limit of Christoffel function.  We also study the relation between the normalized counting measure based on the zeros of the modified average characteristic polynomial and the Christoffel function in limit. The proofs of corresponding theorems with respect to ordinary orthogonal polynomials are based on the three-term recurrence relation. The main point is that exceptional orthogonal polynomials possess at least five-term formulae and so the Christoffel-Darboux formula also fails. It seems that these difficulties can be handled in combinatorial way.
\end{abstract}
\maketitle

\section{Introduction}

Given a (non-trivial) measure $\mu$ supported on a set of the complex plane, say with finite moments, by Gram-Schmidt orthogonalization the sequence of orthogonal polynomials can be generated, $\{P_n\}_{n=0}^\infty \subset L^2(d\mu)$. The kernel of the projection $\pi_n: L^2(d\mu) \to \mathcal{P}_n$, where $\mathcal{P}_n$ is the subspace of polynomials of dimension $n$ is
$$K_n(x,y)=\sum_{k=0}^{n-1}P_k(x)P_k(y).$$
There is an extensive literature of the following two questions investigated in parallel, the asymptotic distribution of zeros of the orthogonal polynomials, $\lim_n\nu_n$, where $\nu_n$ is the normalized counting measure based at the zeros, and the limit of  the weighted reciprocal of the Christoffel function that is the measure
$$d\mu_n(x):=\frac{1}{n}K_n(x,x)d\mu(x).$$
These questions are in close contact with each other and are related to the recurrence relation of orthogonal polynomials.

For instance on the real line if the recurrence coefficients $a_n$ and $b_n$ have limits $a$ and $b$ respectively, the polynomials have asymptotic zero distribution with density $\omega_{a,b}(x)=\frac{1}{\pi\sqrt{(b+2a-x)((x-b+2a)}}$ on the interval $[b-2a,b+2a]$, see \cite{ka} and the references therein.

On the other hand for a $\mu$ supported on the unite circle or on the interval $[-1,1]$ M\'at\'e, Nevai and Totik proved that provided the so-called Szeg\H o condition, $\mu_n$ tends to the equilibrium measures, to the arc- or to the arcsin measure respectively, cf. \cite{mnt}. The asymptotics of $\mu_n$ is investigated in rather general circumstances, we mention here a result due to Totik on arcs and curves, see \cite{t}. Finally it can be shown that the zeros of orthogonal polynomials are the eigenvalues of the operator $\pi_nM_x\pi_n$, where $\pi_n$ is the projection mentioned above and $M_xf=xf$, which ensures that  $(\mu_n-\nu_n)\to 0$ in some sense, see \cite{s}. This observation ensures a direct relation between the two measures derived by ordinary orthogonal polynomials and this technique leads to determinantal point processes.

The aim of the paper is to make small contribution to examination of the problem sketched above in case of exceptional orthogonal polynomials. We focus onto the Christoffel function part.

The notion of exceptional orthogonal polynomials is  motivated by problems in quantum mechanics. It turned up in the last twenty years; one of the earliest papers is \cite{gukm0}, so the topic is fairly new. Owing to the widespread investigation of exceptional systems its literature became rather rich. We mention here just some examples:  \cite{ggm},  \cite{q}, \cite{gukm2}, etc. and see also the references therein.

The study of Christoffel functions of ordinary orthogonal systems depend on the three-term recurrence formula. Exceptional orthogonal polynomials fulfil $2L+1$-term recurrence relations, where $L$ is at least 2. The other important difference is while in ordinary case $xP_n$ can be expressed by a (three-term) linear combination of orthogonal polynomials, in exceptional case instead of $x$ there is a polynomial of degree at least two. This is why most of the methods developed to ordinary orthogonal polynomials are not applicable to the exceptional systems. This investigation based on combinatorial methods inspired by \cite{ha1} and \cite{ha2}. Below we deal with $X1$ Jacobi polynomials in general (generated by one Darboux transform from Jacobi polynomials) and specially $X1$ Jacobi polynomials of codimension 1.

The paper is organized as it follows. In the next section we define exceptional orthogonal polynomials and list up some properties which are important for further investigation. In Section 3 we deal with recurrence relation and give the asymptotics of the recurrence coefficients. In the fourth section we study the asymptotic behavior of the weighted reciprocal of the Christoffel function for codimension 1 exceptional Jacobi case, and in the last section we introduce the modified average characteristic polynomial and study the connection between its zero distribution and $\mu_n$.

\section{Exceptional orthogonal polynomials}

In this section we summarize some properties of exceptional orthogonal polynomials which are important hereinafter. Finally we deal with $X1$-Jacobi polynomials, but our investigation is based on the general characterization of exceptional systems rather than the known examples. The definition is as it follows, cf. \cite[Definition 7.4]{ggm}

Exceptional polynomials means a co-finite real-valued polynomial sequence $\{p_k\}$, where from the sequence of degrees finite many $(k_1, \dots, k_m)$ are missing provided

\noindent (1) The polynomials are eigenfunctions of a differential operator of second order with rational coefficients.

\noindent (2) There is an interval $I$ and a positive weight function $W$ on $I$ with finite moments and at the endpoints of $I$ $pp_kW\to 0$, where $p$ is the coefficient of the second derivative in the differential operator.

\noindent (3) The vector space spanned by the elements of the sequence is dense in the weighted space $L^2(W,I)$.

Second property ensures orthogonality of the exceptional polynomials on $I$ with respect to $W$. Exceptional orthogonal polynomial systems (XPS) are characterized by their construction cf. \cite[Theorem 1.2]{ggm}: all XPS can be obtained by applying finite sequence of Darboux transformations to a classical orthogonal polynomial system.

Below we need the construction thus we summarize the $1$-step Darboux transformation case in short.

Classical orthogonal polynomials $\{P_n^{[0]}\}_{n=0}^\infty$ are eigenfunctions of the second order linear differential operator
$$T[y]=py''+qy'+ry,$$
with eigenvalues $\lambda_n$.
$T$ can be decomposed as
\begin{equation}\label{si}T=BA+\tilde{\lambda}, \ws \mbox{with} \ws A[y]=b(y'-wy), \ws B[y]=\hat{b}(y'-\hat{w}y),\end{equation}
where $b$, $w$ are rational functions and
\begin{equation}\label{s0}\hat{b}=\frac{p}{b}, \ws \ws \hat{w}=-w-\frac{q}{p}+\frac{b'}{b}.\end{equation}
Then the exceptional polynomials are the eigenfunctions of $\hat{T}$ that is the partner operator of $T$ which is
\begin{equation}\label{ka}\hat{T}[y]=(AB+\tilde{\lambda})[y]=py''+\hat{q}y'+\hat{r}y,\end{equation}
where
\begin{equation}\label{hat}\hat{q}=q+p'-2\frac{b'}{b}p, \ws \hat{r}=r+q'+wp'-\frac{b'}{b}(q+p')+\left(2\left(\frac{b'}{b}\right)^2-\frac{b''}{b}+2w'\right)p,\end{equation}
and $w$ fulfils the Ricatti equation
\begin{equation}\label{w}p(w'+w^2)+qw+r=\tilde{\lambda},\end{equation}
cf. \cite[Propositions 3.5 and 3.6]{ggm}.

\eqref{si} and \eqref{ka} ensures that
\begin{equation}\label{sk}\hat{T}AP_n^{[0]}=\lambda_nAP_n^{[0]},\end{equation}
so exceptional polynomials can be obtained from the classical ones by application of (finite many) appropriate first order differential operator(s) to the classical polynomials. This observation motivates the notation below
\begin{equation}\label{A} AP_n^{[0]}=:P_n^{[1]},\end{equation}
(and recursively $A_sP_n^{[s-1]}=:P_n^{[s]}$ in case of $s$ Darboux transforms.)
The degree of $P_n^{[1]}$ is usually greater than $n$, so subsequently $n$ is an index not a degree any more; it shows that $P_n^{[1]}$ is generated from $P_n^{[0]}$.  Since this notation does not follows where the gaps are in the sequence of degrees, it proved to be useful in handling recurrence formulae, cf. \cite{h1}.

$\{P_n^{[1]}\}_{n=0}^\infty$ is an orthogonal system on $I$ with respect to the weight $W:=\frac{pw}{b^2}$, where $w$ is one of the classical weights. Since at the endpoints of $I$ (if there is any) $p$ may have a zero, $b$ can be zero as well here, but it has no zeros inside $I$, otherwise $W$ would not have finite moments. Canceling the zeros at the endpoints we introduce the following notation.

\note

Let $\frac{p}{b}=\frac{\tilde{p}}{\tilde{b}}$,  such that $(\tilde{p},\tilde{b})=1$ and let $\deg\tilde{b}=m$.

\medskip

\remark

\noindent (R1) As $P_n^{[1]}$-s are polynomials for all $n$, applying operator $A$ to $P_0^{[0]}$ it  can be seen that $bw$ must be a polynomial and to $P_1^{[0]}$ shows that $b$ itself is a polynomial too. Thus we can write $w$ as $w=\frac{g}{b}$, where $g$ is also polynomial.

\noindent (R2)
According to \eqref{hat} $P_n^{[1]}$ satisfies the following differential equation:
$$bp y''+\left(b(q+p'-2b'p\right)y'$$ \begin{equation}\label{bhat}+\left(b(r+q')+gp'-b'(q+p')-pb''+2p\frac{b'}{b}(b'-g)\right)y=b(\lambda_n-\tilde{\lambda})y.\end{equation}
Since all the other terms are polynomials $\frac{pb'}{b}(b'-g)P_n^{[1]}$ is also polynomial for all $n$, thus the coefficient of $P_n^{[1]}$ is polynomial. Furthermore denoting by
\begin{equation}\label{L}L[y]:=-2b'py'+\left(gp'-b'(q+p')-pb''+2p\frac{b'}{b}(b'-g)\right)y,\end{equation}
\eqref{bhat} ensures that $LP_n^{[1]}$ is divisible by $\tilde{b}$ for all $n \in \mathbb{N}$.

\noindent (R3) Observe that if $L[f]=Ff'+Gf$ then $L(fg)=fFg'+gLf$.

\medskip

\section{Recurrence relations}

In this section after introducing recurrence relation with respect to exceptional orthogonal polynomials, we investigate the asymptotics of the recurrence coefficients.

Ordinary orthogonal polynomial systems contain polynomials of all degrees which fulfil the three-term recurrence relation
\begin{equation}\label{klrek}xP_n^{[0]}=a_{n+1}P_{n+1}^{[0]}+b_nP_n^{[0]}+a_nP_{n-1}^{[0]}.\end{equation}
Recurrence formulae can be derived for exceptional orthogonal polynomials as well, but these are $2L+1$-term formulae, where $L\ge 2$. As \eqref{klrek} can be rearranged in two ways, for exceptional polynomials there are two different kinds of recurrence relations: with constant- and with variable dependent coefficients. Recalling the Darboux decomposition of the differential operator which generates exceptional polynomials, formulae with dependent coefficients can be derived by operators $A_i$ (cf. \eqref{A}) (see e.g. \cite{stz}, \cite{gugm}, \cite{h1}, \cite{o}, \cite{o2}, \cite{o3}) and with constant coefficient ones by operators $B_i$ (see e.g. \cite{o2}, \cite{o3}, \cite{gukkm}, \cite{o4}, \cite{d}, \cite{mt}, \cite{cggm}). The length of the formula, $L$ in the variable dependent coefficient case depends on the number of Darboux transforms, that is
$$P_n^{[L-1]}=\sum_{k=-L}^Lr_{n,k}^{[L]}(x)P_{n+k}^{[L-1]}.$$
In contrast in the constant coefficient case the length depends on the codimension independently of how many Darboux transforms resulted the codimension of the space of exceptional polynomials. In variable dependent coefficient case it was more convenient to work with monic polynomials, but in constant coefficient case we use the orthonormal form of the polynomials.

Subsequently we deal with exceptional polynomials generated by one-step Darboux transform, so for sake of simplicity further properties are written in this case.
We denote by $P_n^{[i]}$ the orthonormal classical ($i=0$) or exceptional ($i=1$) polynomials. Below we summarize the results of Odake, G´omez-Ullate, Kasman, Kuijlaars, and Milson (see \cite{o2}, \cite{o3}, \cite{gukkm}). Because we need some steps of the proof too, below the sketch of the proofs are also given in brief and for sake of simplicity in one-step Darboux transform case.

Let us introduce the linear space

\begin{equation}\mathcal{P}^{[1]}:=\left\{\sum_{n=0}^Nc_nP_n^{[1]}: c_n\in \mathbb{R}, N\in \mathbb{N}\right\},\end{equation}
and the stabilizer ring
\begin{equation}\mathcal{S}:=\left\{s\in \mathcal{P} : sP_n^{[1]} \in \mathcal{P}^{[1]} \ws \mbox{for all} \ws n\in \mathbb{N}\right\},\end{equation}
where $\mathcal{P}$ denotes the set of polynomials.

\medskip

\noindent {\bf Theorem A.} With the notation above
{Let $s\in \mathcal{P}$. If $s'$ is divisible by $\tilde{b}$, then $s \in \mathcal{S}$.}

\medskip

\proof By (R2) and by comparison of the codimensions of the two sets in question, it can be seen that a polynomial $s \in \mathcal{P}$ belongs to $\mathcal{P}^{[1]}$ if and only if $Ls$ is divisible by $\tilde{b}$. According to (R3) $LsP_n^{[1]}=-2b'ps'P_n^{[1]}+sLP_n^{[1]}$. Thus using (R2) if $s'$ is divisible by $\tilde{b}$ then $LsP_n^{[1]}$ is also divisible by $\tilde{b}$ and so it belongs to $\mathcal{P}^{[1]}$. By orthogonality properties it can be expand as it is given below
\begin{equation}\label{R}QP_n^{[1]}=\sum_{k=-L}^Lu_{n,k}P_{n+k}^{[1]}.\end{equation}

\medskip

\noindent {\bf Corollary B.} {\it If $$BQP_n^{[1]}=\sum_{k=-L}^La_{n,k}P_{n+k}^{[0]},$$ then
$$u_{n,k}=\frac{a_{n,k}}{\lambda_{n+k}-\tilde{\lambda}}, \ws k=-L, \dots , L.$$}

\medskip

\proof $$BQP_n^{[1]}=\sum_{k=-n}^L\frac{a_{n,k}}{\lambda_{n+k}-\tilde{\lambda}}(\lambda_{n+k}-\tilde{\lambda})P_{n+k}^{[0]}=\sum_{k=-n}^L\frac{a_{n,k}}{\lambda_{n+k}-\tilde{\lambda}}BAP_{n+k}^{[0]}$$ $$=\sum_{k=-n}^L\frac{a_{n,k}}{\lambda_{n+k}-\tilde{\lambda}}BP_{n+k}^{[1]}=B\sum_{k=-L}^L\frac{a_{n,k}}{\lambda_{n+k}-\tilde{\lambda}}P_{n+k}^{[1]}.$$
That is
$$0 \equiv B\sum_{k=-L}^L\left(u_{n,k}-\frac{a_{n,k}}{\lambda_{n+k}-\tilde{\lambda}}\right)P_{n+k}^{[1]}=\sum_{k=-L}^L\left(u_{n,k}-\frac{a_{n,k}}{\lambda_{n+k}-\tilde{\lambda}}\right)BAP_{n+k}^{[0]}$$ $$=\sum_{k=-L}^L\left(u_{n,k}-\frac{a_{n,k}}{\lambda_{n+k}-\tilde{\lambda}}\right)(\lambda_{n+k}-\tilde{\lambda})P_{n+k}^{[0]},$$
which ensures that $u_{n,k}-\frac{a_{n,k}}{\lambda_{n+k}-\tilde{\lambda}}=0$ for all $k$.

\medskip

\subsection{Asymptotics of recurrence coefficients}
The purpose of this section is determination of the limits of recurrence coefficients i.e. $\lim_{n\to\infty}u_{n,k}=:U_{|k|}$, cf. \eqref{R}. To this end we return to the expression in Corollary B.

Subsequently we use the simplest element of $\mathcal{S}$:
$$Q(x)=\int^x\tilde{b}, \ws \ws \deg Q=L,$$
where $L=m+1$. Applying (R3), \eqref{s0} and \eqref{si}
$$BQP_n^{[1]}=\hat{b}Q'P_n^{[1]}+QBP_n^{[1]}$$ $$=\tilde{p}P_n^{[1]}+(\lambda_n-\tilde{\lambda})QP_n^{[0]}=\tilde{p}b\left(P_n^{[0]}\right)'-\tilde{p}gP_n^{[0]}+(\lambda_n-\tilde{\lambda})QP_n^{[0]}.$$
Let
\begin{equation}\label{bQ}\tilde{b}(x)=\sum_{k=0}^md_kx^k, \ws \ws Q(x)=\sum_{k=1}^{L}\frac{d_{k-1}}{k}x^k \ws \ws \mbox{and} \ws \ws (\tilde{p}g)(x)=\sum_{k=0}^{m+1}c_kx^k.\end{equation}
In Theorem A the assumption refers to $Q'$ thus we can take $Q$ with an arbitrary constant term. Here the constant is chosen to be zero.

Classical orthonormal polynomials satisfy the following relation (cf. e.g. \cite{sz}).
\begin{equation}\label{kl}p\left(P_n^{[0]}\right)'=A_nP_{n+1}^{[0]}+B_nP_n^{[0]}+C_nP_{n-1}^{[0]},\end{equation}
where $A_n$, $B_n$ and $C_n$ are constants and $p$ is the coefficient of the second derivative in the differential equation. Thus
\begin{equation}\label{BQP}BQP_n^{[1]}=\tilde{b}\left(A_nP_{n+1}^{[0]}+B_nP_n^{[0]}+C_nP_{n-1}^{[0]}\right)-\tilde{p}gP_n^{[0]}+(\lambda_n-\tilde{\lambda})QP_n^{[0]}\end{equation} $$=A_n\sum_{k=0}^md_k\left(x^kP_{n+1}^{[0]} \right)+C_n\sum_{k=0}^md_k\left(x^kP_{n-1}^{[0]} \right)$$ $$+\left(B_nd_0-c_0+\sum_{k=1}^m\left(B_nd_k-c_k+\frac{\lambda_n-\tilde{\lambda}}{k}d_{k-1}\right)x^k+(\lambda_n-\tilde{\lambda})\frac{d_m}{m+1}x^{m+1}\right)P_n^{[0]}.$$
That is knowing the limits of the recurrence coefficients of classical orthonormal polynomials similar limits can be derived to the corresponding exceptional ones.

Wide class of ordinary orthogonal polynomials supported on $[-1,1]$ have recurrence coefficients with limit $\lim_{n\to \infty}a_n=\frac{1}{2}$; $\lim_{n\to \infty}b_n=0$, see \cite[Theorem 4.5.7]{ne}. In classical cases, the situation is easier, after a normalization if it is necessary the limit of the recurrence coefficients; $\lim_{n\to \infty}a_n=:a$; $\lim_{n\to \infty}b_n=:b$ can be got by direct computations. That is \eqref{klrek} can be written as
\begin{equation}xP_n^{[0]}=(a+\varepsilon_{n,1})P_{n+1}^{[0]}+(b+\epsilon_{n,0})P_n^{[0]}+(a+\varepsilon_{n,-1})P_{n-1}^{[0]},\end{equation}
where $\varepsilon_{n,j}$, $\epsilon_{n,j}$ tend to zero. Taking the product of $k$ $a+\varepsilon_{n,j}$ and $b+\epsilon_{n,j}$, where $n-k\le j \le n+k$ the error term can be estimated by $c(k)\varepsilon_{n,0}$, say, where $c(k)$ is a constant depending on $k$. The symmetry of \eqref{klrek} ensures the same symmetry of the iterated recurrence relation. It has the following form:
\begin{equation}\label{itkl}x^kP_n^{[0]}=\sum_{j=-k}^k(s_{k,|j|}+e_{k,n,j})P_{n+j}^{[0]},\end{equation}
where $|e_{k,n,j}|\le c(k)\varepsilon_{n,0}$ that is it tends to zero when $n$ tends to infinity. Substituting it into \eqref{BQP} we have
\begin{equation}\label{BQP2}BQP_n^{[1]}\end{equation} $$=A_n\sum_{k=0}^md_k\sum_{j=-k}^k(s_{k,|j|}+e_{k,n+1,j})P_{n+1+j}^{[0]}+C_n\sum_{k=0}^md_k\sum_{j=-k}^k(s_{k,|j|}+e_{k,n-1,j})P_{n-1+j}^{[0]}$$ $$
+(B_nd_0-c_0)P_{n}^{[0]}+\sum_{k=1}^m\left(B_nd_k-c_k+\frac{\lambda_n-\tilde{\lambda}}{k}d_{k-1}\right)\sum_{j=-k}^k(s_{k,|j|}+e_{k,n,j})P_{n+j}^{[0]}$$ $$+(\lambda_n-\tilde{\lambda})\frac{d_m}{m+1}\sum_{j=-(m+1)}^{m+1}(s_{m+1,|j|}+e_{m+1,n,j})P_{n+j}^{[0]}.$$
In view of \eqref{BQP2} and Corollary B the coefficients in the expansion of $BQP_n^{[1]}$ are
\begin{equation}\label{u1}u_{n,\pm(m+1)}=\frac{1}{\lambda_{n\pm(m+1)}-\tilde{\lambda}}\left((\lambda_n-\tilde{\lambda})\frac{d_m}{m+1}(s_{m+1,|m+1|}+e_{m+1,n,m\pm 1})\right.\end{equation} $$\left. +A_nd_m(s_{m,|m|}+e_{m,n+1,m})\right),$$
and if $1\le |j|\le m$
\begin{equation}\label{u2}u_{n,\pm j}=\frac{1}{\lambda_{n\pm j}-\tilde{\lambda}}\left(A_n\sum_{k=|\pm j-1|}^{m}d_k(s_{k,|j-1|}+e_{k,n+1,j\pm 1})\right.\end{equation}$$\left.+C_n\sum_{k=|\pm j+1|}^{m}d_k(s_{k,|j+1|}+e_{k,n-1,j\pm 1})\right.$$ $$\left.+\sum_{k=\max\{1,|j|\}}^{m}\left(B_nd_k-c_k+(\lambda_n-\tilde{\lambda})\frac{d_{k-1}}{k}\right)(s_{k,|j|}+e_{k,n,j})\right.$$ $$\left.+(\lambda_n-\tilde{\lambda})\frac{d_m}{m+1}(s_{m+1,|j|}+e_{m+1,n,j})\right).$$

\noindent To finish this general part we compute $s_{k,j}$.

\medskip

\begin{lemma}\label{sj} In \eqref{itkl} the iterated coefficients in limit are
\begin{equation}\label{s} s_{k,j}=\sum_{i=0}^{\left[\frac{k-|j|}{2}\right]}\binom{k}{|j|+2i}\binom{|j|+2i}{i}a^{|j|+2i}b^{k-|j|-2i}.\end{equation}\end{lemma}

\medskip

\remark
Certainly similar results can be found in the literature, see e.g. \cite{bc} with certain Freud weights. We introduce here the method we use mainly in the next section to get the result in the required form.

\medskip

\proof
Recalling that $s_{k,j}$ is the main part of the coefficient of $P_{n+j}^{[0]}$ in the expansion of $x^kP_{n}^{[0]}(x)$ we can compute it on an oriented weighted graph $G=(V,E)$, where the vertices $V=\mathbb{N}^2$ and the edges
$$E:=\left\{e:(n,m)\to (n+1,m+k), n,m\in \mathbb{N}, \ws -1\le k \le 1\right\}.$$
The weight on an edge is given by
$$w(e)=\left\{\begin{array}{ll}a, \ws \mbox{if} \ws k=1 \ws \mbox{or} \ws k=-1\\ b, \ws \mbox{if} \ws k=0.\end{array}\right.$$
We start from zero level, say  and we have to reach level $j$ (or $-j$) in $k$ steps. As we can step forward (on the same level) with weight $b$, up and down with weight $a$, on the $k$-long path we can take $j+2i$ steps up or down which contains $j+i$ steps upwards and $i$ downwards, and the remainder ones forwards. (For $-j$ we get the same.) That is
$$ s_{k,j}=\sum_{\gamma: (0,n) \to (k,n+j)}\prod_{e\in\gamma} w(e)=\sum_{i=0}^{\left[\frac{k-|j|}{2}\right]}\binom{k}{|j|+2i}\binom{|j|+2i}{i}a^{|j|+2i}b^{k-|j|-2i},$$
which is just the statement.

\noindent In the rest of this section we compute the asymptotic recurrence coefficients $U_{|k|}$.

\subsection{ X1 Jacobi polynomials}
Let the $n^{th}$ Jacobi polynomial defined by Rodrigues' formula:
$$(1-x)^{\alpha}(1+x)^{\beta}P_n^{\alpha,\beta}(x)=\frac{(-1)^n}{2^nn!}\left((1-x)^{\alpha+n}(1+x)^{\beta+n}\right)^{(n)},$$
where $\alpha, \beta >-1$.
$$p_n^{\alpha,\beta}=\frac{P_n^{\alpha,\beta}}{\sqrt{\sigma_n}},$$
the orthonormal Jacobi polynomials which fulfil the following differential equation (cf. \cite[(4.2.1)]{sz})
\begin{equation}\label{jd}(1-x^2)y''+(\beta-\alpha-(\alpha+\beta+2)x)y'+n(n+\alpha+\beta+1)y=0.\end{equation}
With
\begin{equation}\label{jac}P_{n}^{[0]}=p_n^{\alpha,\beta},\end{equation}
so taking into consideration that $\frac{\sigma_{n-1}}{\sigma_n}\sim 1$, etc. (for $\sigma_n$ see \cite[(4.3.4)]{sz}), according to \cite[(4.5.5)]{sz}
\begin{equation}\label{ABC} A_n\sim\frac{1}{2}n, \ws \ws B_n \sim \frac{\alpha-\beta}{2} \ws \ws C_n\sim -\frac{1}{2}n.\end{equation}
$A_n:=A_n(\alpha,\beta)$, etc. - for sake of simplicity we omit $\alpha$ and $\beta$ and $X \sim Y$ means that $\frac{X}{Y}$ is between two positive constants.
Again after normalization in view of \cite[(4.5.1)]{sz}
\begin{equation}\label{ab}a_n\to \frac{1}{2}, \ws \ws b_n\to 0.\end{equation}

\medskip

Applying Lemma \ref{sj} with the values of \eqref{ab} we immediately have

\begin{lemma}\label{sjj}

$$s_{k,j}=\left\{\begin{array}{ll}\binom{k}{\frac{k-j}{2}}\frac{1}{2^k}, \ws \mbox{if $k-j$ is even}\\0, \ws \mbox{if $k-j$ is odd.}\end{array}\right.$$
\end{lemma}

\medskip

Now we are in position to compute the limit coefficients. Since $\tilde{\lambda}$ does not depend on $n$, by \eqref{jd} $\lambda_{n\pm j}-\tilde{\lambda}\sim n^2$. Taking into consideration this and substituting \eqref{ab}, \eqref{ABC} and Lemma \ref{sjj} into \eqref{u1} and \eqref{u2} we arrive to

\begin{proposition}\label{p1}
If $L=m+1$ is even or odd respectively
\begin{equation}\label{u3}U_{|j|}=\left\{\begin{array}{ll}\sum_{p=\max\{l,1\}}^{\frac{m+1}{2}}\frac{d_{2p-1}}{2p}\binom{2p}{p-l}\frac{1}{2^{2p}}, \ws\ws \mbox{if}\ws \ws |j|=2l\\
\sum_{p=l}^{\frac{m}{2}}\frac{d_{2p}}{2p+1}\binom{2p+1}{p-l}\frac{1}{2^{2p+1}}, \ws\ws \mbox{if}\ws \ws |j|=2l+1. \end{array}\right.
\end{equation}\end{proposition}

\medskip

\section{$X1$ Jacobi polynomials of codimension 1}

For ordinary orthogonal polynomials it is proved in rather general circumstances that $\frac{1}{n}K_n(x,x)w(x)dx$ tends to the equilibrium measure of the (essential) support of $\mu$ in some sense, where $w$ is the Radon-Nikodym derivative of $\mu$ and  $K_n(x,x)$ the Christoffel-Darboux kernel, see e.g. \cite{s} and the references therein. All he proofs of these types of theorems use the three-term recurrence relation of orthogonal polynomials. In this section we show in the simplest, one codimensional case that is with five-term recurrence formula, that the statement above still fulfils.

With the notation above let us choose
$$\tilde{b}(x)=d_1x+d_0, \ws \mbox{that is} \ws \ws  Q(x)=\frac{d_1}{2}x^2+d_0x.$$
Here we choose the constant term of $Q$ to $0$ again. To define $\tilde{b}$ monic polynomials are also appropriate. Since in the known examples $\tilde{b}$ is not monic (cf. e.g. \cite{gumm}) we take it in the form given above and it does not cause any further difficulty. Let
$$K_n(x,y):=\sum_{k=0}^{n-1}P_k^{[1]}(x)P_k^{[1]}(y),$$
where $P_k^{[1]}(x)$ are the $k^{th}$ exceptional Jacobi polynomials orthonormal with respect to $W$ on $[-1,1]$. We introduce the following notation:
\begin{equation}\label{mn}d\mu_N(x)=\frac{1}{N}K_N(x,x)W(x)dx,\end{equation}
$$\langle f,g \rangle :=\int_{-1}^1f(x)g(x)W(x).$$
The equilibrium measure of $[-1,1]$, $\mu_e$ is defined by
$$d\mu_e(x)=\omega(x)dx=\frac{1}{\pi}\frac{1}{\sqrt{1-x^2}}dx.$$
Let us recall the Wallis' formula:
$$\int_{-1}^1x^ld\omega(x)=\left\{\begin{array}{ll}\binom{2k}{k}\frac{1}{2^{2k}}, \ws l=2k,\\0, \ws l=2k+1.\end{array}\right.$$

\medskip

\begin{theorem}\label{Tw}  If $\mu_N$ is defined as in \eqref{mn} by $X_1$-Jacobi polynomials of codimension 1, on $[-1,1]$
$$\lim_{N\to\infty}\mu_N=\mu_e$$
in weak-star sense.
\end{theorem}

\medskip

To prove Theorem \ref{Tw} we need the following Lemma.

\medskip

\begin{lemma}\label{le2} For all $k=0,1, \dots$
\begin{equation}\label{wk}\lim_{n \to \infty}\langle Q^kP_n^{[1]},P_n^{[1]}\rangle=\int_{-1}^1Q^k(x)\omega(x)dx.\end{equation}\end{lemma}

\medskip

\proof
First we compute the effect of $\omega$ to $Q^k$, $k=0,1,\dots$. Binomial theorem and Wallis' formula give
\begin{equation}\label{QQ}\int_{-1}^1Q^k(x)\omega(x)dx\end{equation} $$=\sum_{j=0}^k\binom{k}{j}d_0^jd_1^{k-j}\frac{1}{2^{k-j}}\int_{-1}^1x^{2k-j}(x)d\omega(x)=\sum_{i=0}^{\left[\frac{k}{2}\right]}\binom{k}{2i}\binom{2(k-i)}{k-i}\frac{d_0^{2i}d_1^{k-2i}}{2^{3k-4i}}.$$

On the other hand, according to Theorem A for all (fixed) $k$  there exist coefficients $c_{n,j}^{(k)}$ such that
$$Q^kP_n^{[1]}=\sum_{j=-2k}^{2k}c_{n,j}^{(k)}P_{n+j}^{[1]},$$
thus
$$\lim_{n \to \infty}\langle Q^kP_n^{[1]},P_n^{[1]}\rangle=\lim_{n \to \infty}c_{n,0}^{(k)}=:c^{(k)},$$
where the limit exists with regard to Proposition \ref{p1}.
We can proceed as previously, that is recalling the oriented weighted graph $G=(V,E)$ with $V=\mathbb{N}^2$, we define edges as
$$E:=\left\{e:(n,m)\to (n+1,m+k), n,m\in \mathbb{N}, \ws -2\le k \le 2\right\},$$
and weights on the edges as
$$w(e)=\left\{\begin{array}{ll}U_0,\ws \mbox{if} \ws k=0\\ U_1, \ws \mbox{if} \ws k=1 \ws \mbox{or} \ws k=-1\\ U_2, \ws \mbox{if}  \ws k=2 \ws \mbox{or} \ws k=-2 ,\end{array}\right.$$
where, according to Proposition \ref{p1}
\begin{equation}\label{U} U_0=\frac{d_1}{4}, \ws \ws U_1=\frac{d_0}{2}, \ws \ws U_2=\frac{d_1}{8}.\end{equation}
Thus $c^{(k)}$ can be expressed as the collection of all possible weighted $k$-long paths from $P_n^{[1]}$ to $P_n^{[1]}$,
$$c^{(k)}=\sum_{\gamma(0,n)\to (k,n)}\prod_{e\in\gamma}w(e).$$
That is in each step we can move forwards, denoted by $(0)$ (with weight $U_0$); one level up or down denoted by $(\pm 1)$ (with weight $U_1$) and two levels up or down denoted by $(\pm 2)$ (with weight $U_2$). Since finally we have to arrive to the starting level it is clear that we need even number of $(\pm 1)$. Let the number of $(\pm 1)$ is $2i$. Then we can choose $s+m$ pieces of $(2)$ (that is two levels big steps upwards)  and $m$ $(-2)$. The rest of the steps are zeros (i.e. forwards). That is
\begin{equation}\label{ck}c^{(k)}=\sum_{i=0}^{\left[\frac{k}{2}\right]}\binom{k}{2i}\end{equation} $$\times\sum_{s=0}^{\min\{i,k-2i\}}\binom{2i}{s+i}\sum_{m=0}^{\left[\frac{k-2i-s}{2}\right]}\binom{k-2i}{s+2m}\binom{s+2m}{m}\frac{d_0^{2i}d_1^{k-2i}}{2^{2k-2i+2m+ \max\{0,s-1\}}}.$$
Indeed we can choose the $2i$ $(\pm 1)$-s in $\binom{k}{2i}$ ways, and to neutralize the $s$ pieces of $(2)$ we need $2s$ $(-1)$ and we have further $\frac{2i-2s}{2}$ $(-1)$, that is the number of the choices of $(-1)$-s is $\binom{2i}{s+i}$. Then we choose from the remainder $k-2i$ elements the $s+2m$ elements with absolute value $2$, and from these the $m$ pieces of $(-2)$. If $s>0$, we can do the same with $s$ $(-2)$-s, etc. and it gives a factor $2$ in these cases. Now taking into consideration the weights we arrive to the formula above.

Notice that \eqref{ck} can be rearranged as
$$c^{(k)}=\sum_{i=0}^{\left[\frac{k}{2}\right]}\binom{k}{2i}\frac{d_0^{2i}d_1^{k-2i}}{2^{2k-2i}}S_{k,i},$$
where
$$S_{k,i}=\sum_{s=0}^{\min\{i,k-2i\}}\binom{2i}{s+i}\sum_{m=0}^{\left[\frac{k-2i-s}{2}\right]}\binom{k-2i}{s+2m}\binom{s+2m}{m}\frac{1}{2^{2m+ \max\{0,s-1\}}}.$$
In order to prove that $c^{(k)}=\int_{-1}^1Q^k(x)d\omega(x)$, it is enough to show that for all $0\le i \le \left[\frac{k}{2}\right]$,
\begin{equation}\label{S} S_{k,i}=\binom{2(k-i)}{k-i}\frac{1}{2^{k-2i}}.\end{equation}

First we observe that if $i=0$ and $k$ is arbitrary, then
$$S_{k,0}=\sum_{m=0}^{\left[\frac{k}{2}\right]}\binom{k}{2m}\binom{2m}{m}\frac{1}{2^{2m}}$$ \begin{equation}\label{S0}=\int_{-1}^1(1+x)^kd\omega(x)=\frac{2^{k+1}}{\pi}\int_0^{\frac{\pi}{2}}\cos^{2k}\xi d\xi=\binom{2k}{k}\frac{1}{2^{k}}.\end{equation}
If $k$ is even and $i=\frac{k}{2}$, then on both of the sides of  \eqref{S} stands $\binom{k}{\frac{k}{2}}$.

To show \eqref{S} we use induction from $(k,i)$ to $(k+1, i+1)$. Notice that replacing $(k+1, i+1)$ instead of $(k,i)$ the right-hand side of \eqref{S} doubles up. We show the same for the left-hand side of \eqref{S}.

We can consider $S_{k,i}$ as all the weighted paths from zero level to zero level in $k$ steps with fixed $2i$ $(\pm 1)$-s and with the following new weights: $(\pm 2)$ have weight $\frac{1}{2}$ and $(0)$, $(\pm 1)$ have weight $1$. In $S_{k+1,i+1}$ we have one more step and two more $(\pm 1)$-s. Thus from a paths of $S_{k,i}$ we can get a path of $S_{k+1,i+1}$ if we replace a $(2)$ by two $(-1)$-s (or a $(-2)$ by two $(1)$-s), or a $(0)$ by $(1)$ and $(-1)$ and there are no other ways. Observe that all these steps double up the weighted sum, because if we replace a $(2)$ by $(-1)$-s  we replace $\frac{1}{2}$ by $1$ and when we replace a $(0)$, we can change the places of $(1)$ and $(-1)$, so we get two different paths. Proceeding this way and then omitting the paths we have got more than one times, we get all the paths of $S_{k+1,i+1}$. Conversely, taking two paths from $S_{k+1,i+1}$ such that one can be derived from the other by a change of a pair of $(\pm 1)$ and replacing them by a $(0)$, or two $(1)$-s by a $(-2)$, etc., omitting the repetitions again we get all the paths of $S_{k,i}$, thus $S_{k+1,i+1}=2S_{k,i}$ and together with \eqref{S0} it gives the statement.

\medskip

\proof (of Theorem \ref{Tw}) By arithmetic mean law Lemma \ref{le2} ensures for all $k=0,1, \dots$
$$\lim_{N\to\infty}\int_{-1}^1Q^k(x)d\mu_N(x)=\int_{-1}^1Q^k(x)d\omega(x).$$
The linear space
$$A:=\mathrm{span}\left\{Q^k :k\in\mathbb{N}\right\}$$
is also an algebra with unit element. Recalling that $W$ has finite moments, and $\tilde{b}^2$ is in the denominator of $W$, the zero of $\tilde{b}$ is out of the interval of orthogonality, that is $\tilde{b}=Q'$ has constant sign here and so $Q$ is monotone on $[-1,1]$ which ensures that $A$ separates the points of $[-1,1]$. So by Stone-Weierstrass theorem $A$ is dense in $C[-1,1]$. As $\mu_N$ for all $N$ and $\omega$ are probability measures, the theorem follows from density.

\section{"Average characteristic polynomial"}

In this section we look at $\mu_N$ (see \eqref{mn}) from another perspective. It shows up in connection with determinantal point processes. Without going into this theory deep, we sketch how it relates to average characteristic polynomial, $\chi_N$.

In ordinary case this correspondence ensures results about limit zero distribution of orthogonal polynomials. Although the known examples show that normalized counting measures based on the zeros of exceptional polynomials tend to the equilibrium measure too (see e.g. \cite{gumm}, \cite{h}, \cite{km}, \cite{h2}, \cite{l}), the situation in general case is a little bit more complicated.

To get relationship between $\mu_N$ and $\chi_N$ similar to the ordinary case, we have to modify the average characteristic polynomial - we denote the modified version with $\chi_N$ again.

Let $\{P_n\}_{n=0}^{\infty}$ be an ordinary or exceptional orthonormal system of polynomials on a real interval $I$ with respect to an appropriate measure $\mu$ or weight function $W$ (which is the Radon-Nikodym derivative of $\mu$), that is orthonormal in $L^2(\mu)$ and $\deg P_n\ge n$. Let
$$\mathcal{P}_N:=\mathrm{span}\left\{P_n: n=0, \dots , N-1\right\}, \ws \ws \mathcal{P}=\bigcup_{N\ge 1}\mathcal{P}_N$$
and as above
$$K_N(x,y):=\sum_{k=0}^{N-1}P_k(x)P_k(y).$$
We define a polynomial $Q$ and $u_{n,j}$ such that
\begin{equation}\label{u9}QP_n \in \mathcal{P}\ws \ws \forall \ws n\in\mathbb{N} \ws \ws \mbox{and}\ws \ws QP_n=\sum_{j=-L}^Lu_{n,j}P_{n+j}\end{equation}
with some fixed $L$ depending only on $Q$. We also define two mappings acting on $L^2(\mu)$
$$\pi_N: f(x)\mapsto \int_IK_N(x,t)f(t)d\mu(t)\in \mathcal{P}_N, \ws \ws \ws M: f(x)\mapsto Q(x)f(x).$$
Note that in ordinary case $\deg P_n=n$, $Q(x)=x$ as in \cite{ha1} and \cite{ha2}.

$x_1, \dots, x_N$ are random variables, the joint probability distribution on $\mathbb{R}^N$ is
\begin{equation}\label{ro}\varrho_N(x_1, \dots ,x_N)=c(N)\det|K_n(x_i,x_j)|_{i,j=1}^N\prod_{i=1}^NW(x_i),\end{equation}
where $c(N)$ is a normalization factor. The expectation $\mathbb{E}$ refers to $\varrho$.

Let $f$ be real valued, bounded and measurable function with bounded support contained by $I$ and assume that
\begin{equation}\label{jo}\sum_{N=0}^\infty\frac{\|f\|^N_{\infty}}{N!}\int_{I^N}\varrho_N(x_1, \dots ,x_N)dx_1\dots dx_N<\infty.\end{equation}
Then
\begin{equation}\label{E}\mathbb{E}\left(\sum_{i_1\neq \dots \neq i_k}f(x_{i_1})\dots f(x_{i_k})\right)\end{equation}$$=\int_{I^k}f(x_1) \dots  f(x_k)\det|K_n(x_i,x_j)|_{i,j=1}^k\prod_{i=1}^kd\mu(x_i),$$
see e.g. \cite[Proposition 2.2]{j}.

The empirical distribution
$$\hat{\mu}_N=\frac{1}{N}\sum_{i=1}^N\delta_{x_i}.$$
Specially if $k=1$ $n=N$ \eqref{E} becomes
\begin{equation}\label{E1}\mathbb{E}\left(\int f\hat{\mu}_N\right)=\int_I f(x) \frac{1}{N}K_N(x,x)W(x)dx=\int_Ifd\mu_N.\end{equation}
We also need $\hat{\mu}_N^Q=\frac{1}{N}\sum_{i=1}^N\delta_{Q(x_i)}$ and the modified average characteristic polynomial
$$\chi_N(z):=\chi_N(z)^Q(z)=\mathbb{E}\left(\prod_{i=1}^N(z-Q(x_i))\right).$$
Denote by $z_i$ the zeros of $\chi_N(z)$ and define
$$\nu=\frac{1}{N}\sum_{i=1}^N\delta_{z_i}.$$

After introducing notation above we are in position to state the main result of this section.

\medskip

\begin{theorem}\label{T2}
If there is a $B$ such that $|u_{n,j}|\le B$ for all $n,j$, then for all $l\ge 0$
$$\lim_{N\to \infty}\left|\mathbb{E}\left(\int x^ld\hat{\mu}_N^Q(x)\right)- \int x^ld\nu_N(x)\right|=0.$$\end{theorem}

 \medskip

 Let $z_i=Q(y_i)$, $i=1, \dots , N$, $\tilde{\nu}_N=\frac{1}{N}\sum_{i=1}^N\delta_{y_i}$, that is $\int Q^l(y)d\tilde{\nu}_N(y)=\int z^ld\nu_N(z)$. Proposition \ref{p1} and Theorem \ref{T2} immediately ensures the next corollary.

 \medskip

\begin{corollary}  In the $X_1$ Jacobi case for all $l\ge 0$
$$\lim_{N\to \infty}\left|\int_{-1}^1Q^ld\mu_N- \int Q^ld\tilde{\nu}_N\right|=0.$$\end{corollary}

\medskip

\begin{lemma}\label{tr1}

\begin{equation}\label{l1} \mathbb{E}\left(\int x^ld\hat{\mu}_N^Q(x)\right)=\frac{1}{N}\mathrm{Tr}(\pi_NM^l\pi_N),\end{equation}
\begin{equation}\label{l2}\int x^ld\nu_N(x)=\frac{1}{N}\mathrm{Tr}\left((\pi_NM\pi_N)^l\right).\end{equation} \end{lemma}

\medskip

\proof

To apply \eqref{E} first we have to check the conditions assumed above in our case. We take $f=Q|_{[-1,1]}$, which has the required properties. Orthogonality together with Dyson's theorem (see \cite[Theorem 5.14]{meh}) ensures the convergence in \eqref{jo}.

\noindent \eqref{l1}: In view of \eqref{E1}
$$\mathbb{E}\left(\frac{1}{N}\sum_{i=1}^NQ^l(x_i)\right)=\frac{1}{N}\sum_{i=1}^N\int_IQ^l(x)\frac{1}{N}K_N(x,x)W(x)dx$$ $$=\frac{1}{N}\sum_{k=0}^{N-1}\int_IQ^l(x)P_k^2(x)W(x)dx=\frac{1}{N}\sum_{k=0}^{N-1}\langle \pi_NM^l\pi_NP_k,P_k\rangle=\frac{1}{N}\mathrm{Tr}(\pi_NM^l\pi_N).$$

\noindent \eqref{l2}: Applying \eqref{E}
$$\mathbb{E}\left(\prod_{i=1}^N(z-Q(x_i)\right)=z^N+\sum_{k=1}^N\frac{(-1)^kz^{N-k}}{k!}\mathbb{E}\left(\sum_{i_1\neq \dots \neq i_k}Q(x_{i_1}) \dots Q(x_{i_k})\right)$$ $$=z^N+\sum_{k=1}^N\frac{(-1)^kz^{N-k}}{k!}\int_{I^k}Q(x_1)\dots Q(x_k)\det|K_N(x_i,x_j)|_{i,j=1}^k\prod_{i=1}^kW(x_i)dx_i.$$
Considering the Fredholm's expansion of integral operator $\pi_NM\pi_N$ acting on $\mathcal{P}_N$ with kernel $Q(y)K_N(x,y)$ we have
$$\det (z-\pi_NM\pi_N)=z^N+\sum_{k=1}^N\frac{(-1)^kz^{N-k}}{k!}\int_{I^k}\det|Q(x_j)K_N(x_i,x_j)|_{i,j=1}^k\prod_{i=1}^kW(x_i)dx_i,$$
cf. \cite{ggk}, from which follows the next equality
$$\det (z-\pi_NM\pi_N)=\chi_N(z):=\mathbb{E}\left(\prod_{i=1}^N(z-Q(x_i)\right).$$
That is the roots of $\chi_N(z)$ are the eigenvalues of $\pi_NM\pi_N$, so $\frac{1}{N}\mathrm{Tr}\left((\pi_NM\pi_N)^l\right)=\frac{1}{N}\sum_{i=1}^Nz_i^l=\int x^ld\nu_N(x)$.

 \medskip

Let
$$D_N:=\{(m,n)\in \mathbb{N}^2: n\ge N\}.$$

The proof of the next lemma is independent of the form of $Q$. The weight of a path, $w(\gamma)=\prod_{e\in\gamma}w(e)$.

 \medskip

\begin{lemma}\label{tr2}\cite[Lemma 2.5, Lemma 2.6]{ha1} For all $l\in\mathbb{N}$
\begin{equation}\label{w1} \mathrm{Tr}(\pi_NM^l\pi_N)=\sum_{k=0}^{N-1}\sum_{\gamma:(0,k) \to (l,k)}w(\gamma),\end{equation}
\begin{equation}\label{w2}\mathrm{Tr}\left((\pi_NM\pi_N)^l\right)=\sum_{k=0}^{N-1}\sum_{\gamma:(0,k) \to (l,k) \atop \gamma\cap D_N=\emptyset}w(\gamma).\end{equation} \end{lemma}

 \medskip

 To make the discussion self-contained we insert the proof of Theorem \ref{T2} in brief, cf. \cite{ha2}.

\proof (of Theorem)
In view of Lemma \ref{tr1} and \ref{tr2} we estimate
\begin{equation}\label{sz}\frac{1}{N}\sum_{k=0}^{N-1}\sum_{\gamma:(0,k) \to (l,k) \atop \gamma\cap D_N\neq \emptyset}w(\gamma).\end{equation}
According to \eqref{u9}, if $k<N-lL$ then $\gamma\cap D_N=\emptyset$. By the assumption on $u_{k,j}$ the weight of an $l$-long path is at most $B^l$ and the vertices of admissible paths are contained by the set $\{(n,m): 0\le n\le l, \ws N-l\le m\le N+l\}$ which has at most $(2l)^l$ elements.
That is
$$\frac{1}{N}\sum_{k=0}^{N-1}\sum_{\gamma:(0,k) \to (l,k) \atop \gamma\cap D_N\neq \emptyset}w(\gamma)\le \frac{(2lB)^l}{N}.$$

\medskip

\remark
By Theorem \ref{Tw} and Theorem \ref{T2} with an appropriately defined $Q$  in $X_1$ Jacobi case of codimension 1 the following relation can be derived
$$\tilde{\nu}_N|_{[-1,1]} \xrightarrow[]{w^*}\mu_e.$$
Indeed, Theorem \ref{T2} ensures that there is a compact interval $[a,b]$ which contains $N-o(N)$ $z_i$ eigenvalues for all $N$. Define $Q^*:=AQ+B$ such that $Q^*([-1,1])=[a,b]$. Thus the main part of $\tilde{\nu}_N$ is supported on $[-1,1]$. Applying the binomial theorem to both of the sides of Lemma \ref{le2} we get the same result with $Q^*$ as well. So the density argument mentioned above leads to the weak star convergence of the restriction of $\tilde{\nu}_N$.

\medskip

\noindent \small{Department of Analysis, \newline
Budapest University of Technology and Economics}\newline
\small{g.horvath.agota@renyi.mta.hu}

\end{document}